\documentclass[10pt]{article}

\usepackage{endfloat}
\usepackage{amssymb}
\usepackage[centertags]{amsmath}
\usepackage{amsthm,bm}
\usepackage{cancel}  %
\usepackage{mathdots} %
\usepackage{eucal}   
\usepackage{mathrsfs} 
\usepackage[numbers,comma]{natbib}
\usepackage{comment, url, graphicx, graphics, bm}
\usepackage{hyperref}
\usepackage{changes}
\def\o*{o_{{\!}_{P^*}}}
\def\O*{\cal O_{{\!}_{P^*}}}

\def\.{\mbox{.}}

\def\le{\leqslant}
\def\ge{\geqslant}
\def\sfrac(#1,#2){\mbox{$\frac{#1}{#2}$}}

\def\tr{{\mbox{\tiny{$\mathrm{T}$}}}}

\def\bernstein{Bernstein }

\newtheorem{theorem}{Theorem}

\newtheorem{lemma}[theorem]{Lemma}

\newtheorem{rem}{Remark}[section]
\title{Improvement of a Theorem of Lorentz (1963) and its Generalization to the Multivariate Case}

\author{Zhong Guan\\Department of Mathematical Sciences\\ Indiana University South Bend\\ South Bend, IN 46634, USA\\
zguan@iusb.edu\\
Tao Wang\\
School of Mathematical Sciences\\ Harbin Normal University\\Harbin, 150025, China}

\begin{document}
\maketitle

\section*{Abstract} In this short note we have proved an enhanced version of a theorem of Lorentz \cite{Lorentz-1963-Math-Annalen} and its generalization to the multivariate case which gives a non-uniform estimate of degree of approximation by a polynomial with positive coefficients.
The performance of the approximation  at the vertices of $[0,1]^d$ is more precisely characterized by the improved result and its multivariate generalization.
The latter provides mathematical foundation on which multivariate density approximation by a polynomial with positive coefficients can be established.

\noindent\textbf{Key Words and Phrases:} Degree of Approximation, Non-uniform estimate, Polynomial with positive coefficients.\\
2010 MSC: 41A10


\section{Introduction}
The {\em polynomial of degree $n$ with positive coefficients} studied by Lorentz \cite{Lorentz-1963-Math-Annalen} can be uniquely represented as $P_n(x)=\sum_{k=0}^n b_k p_{nk}(x)$, $b_k\ge 0$,
where $p_{nk}(x)={n\choose k}x^k(1-x)^{n-k}$, $k=0,\ldots,n$, which are linearly independent.  An example is the {\em Bernstein polynomial} \citep{Bernstein},
  $B_n^f(x)=\sum_{k=0}^n f(k/n) p_{nk}(x)$, for a nonnegative function $f$ defined on $[0,1]$.
Recently, Bernstein polynomial and its generalization, polynomial with positive coefficients, find their important applications in nonparametric statistics  for estimating an unknown probability density function $f$ which is assumed to have support on a closed interval \citep{Guan-jns-2015}. The traditional  density estimate is the
kernel density  \citep{Rosenblatt-1956-ann-math-stat,Parzen-1962} which  targets the convolution of the unknown density function and a scaled kernel function,  instead of the unknown density $f$,   and suffers from serious boundary effects in the sense that the estimation errors are usually large at the endpoints of the  {\em closed} supporting interval. In order that the method of Guan \cite{Guan-jns-2015}, which is much better than the kernel method, can be applied to estimate the multivariate density functions, we shall extend results of Lorentz \cite{Lorentz-1963-Math-Annalen} to the multivariate case in this note. Even better, we actually are able to prove an improved generalization of Lorentz's result \cite{Lorentz-1963-Math-Annalen} which clarifies mathematically why the method  of \citep{Guan-jns-2015} and other methods \citep{Vitale1975,Guan:Wu:Zhao-2008} using Bernstein polynomials can very much reduce the boundary effects.
\section{Multivariate Polynomials
with Positive Coefficients}
Throughout the paper, we use bold face letters to denote vectors. For example, $\bm x = (x_1, \ldots, x_d)^\tr$ is a $d$-dimensional vector. We denote the taxicab norm  by $|\bm x|=\sum_{i=1}^d|x_i|$.
 Inequality $\bm x \le \bm y$ is understood componentwise, i.e., $x_j \le y_j$ for all $j = 1, \ldots , d$. The strict inequality
$\bm x < \bm y$ means  $\bm x \le \bm y$  but $\bm x \ne \bm y$.
Let  $C^{(k)}[0,1]^d$ denote the class of functions $f$ on the unit hypercube $[0,1]^d$ that have continuous partial derivatives
$f^{(\bm l)}(\bm t)\equiv\partial^{\langle\bm l\rangle} f(\bm t)/\partial t_1^{l_1}\cdots\partial t_d^{l_d}$, where $0\le \langle\bm l\rangle \equiv l_1+\cdots+l_d \le k$ and $C[0,1]^d=C^{(0)}[0,1]^d$.
Define the modulus of continuity for derivative $f^{(\bm r)}$ by $\omega_{\bm r}(h)=\omega(f^{(\bm r)},h)$,
where  $\omega(f,h)=\max_{|\bm s-\bm t|<h}|f(\bm s)-f(\bm t)|$, $h>0.$
For each $r\ge 0$, let $\omega^{(r)}(h)=\max_{\langle \bm r \rangle=r}\omega_{\bm r}(h)$.

 One version of the {\em multivariate \bernstein polynomial} approximation \citep{Hildebrandt-Schoenberg-1933,Butzer-2d-Bernstein-poly-1953} for $f\in C[0,1]^d$ is defined as
 $ B_{\bm n}^f(\bm x) = \sum_{\bm i=0}^{\bm n}
f (\sfrac(\bm i, \bm n) ) \cdot\prod_{j=1}^d p_{n_ji_j}(x_j)= \sum_{\bm i=0}^{\bm n}
f(\sfrac(\bm i, \bm n) ) \cdot p_{\bm n,\bm i}(\bm x)$,
where $\bm x=(x_1,\ldots,x_d)\in[0,1]^d$, $\bm n=(n_1,\ldots,n_d)$, $\bm i=(i_1,\ldots,i_d)$, 
$\sfrac(\bm i, \bm n)=(\sfrac(i_1,n_1),\ldots,\sfrac(i_d,n_d))$, $p_{\bm n,\bm i}(\bm x)=\prod_{j=1}^d p_{n_ji_j}(x_j)$,
 and 
  $\sum_{\bm i=0}^{\bm n}=\sum_{i_1=0}^{n_1}\,\cdots   \sum_{i_d=0}^{n_d}$.
Note that $\sum_{\bm i=0}^{\bm n} p_{\bm n,\bm i}(\bm x)= \prod_{j=1}^d\sum_{\bm i=0}^{\bm n} p_{n_ji_j}(x_j)=1.$
The {\em multivariate polynomial with positive coefficients} can  then  be defined as
$P_{\bm n}(\bm x) = \sum_{\bm i=0}^{\bm n}
a(\bm i) \cdot p_{\bm n,\bm i}(\bm x)
$ with $a(\bm i)= a(i_1,\ldots,i_d)\ge 0$.
\section{Degree of Approximation}

To generalize the results of  Lorentz  \cite{Lorentz-1963-Math-Annalen} in a convenient way, we define
 $\Lambda_{r}=\Lambda_{r}^{(d)}(m, \mathscr M_{r})$, where $\mathscr M_0=\mathscr M_1=M_0$ and
$\mathscr{M}_{r}=(M_0, M_{\bm i},  2\le\langle\bm i\rangle \le  r)$ for $r\ge 2$,  as the class of
 functions $f(\bm t)$ in $C^{(\bm r)}[0,1]^d$ such that
    $m\le f(\bm t)\le M_0$, $|f^{(\bm i)}(\bm t)| \le  M_{\bm i}$, $\bm t\in [0,1]^d$,
for some $m>0$, $M_{\bm i}$,   and all $\bm i$ such that  $2\le \langle\bm i\rangle \le r$.

Using the  notations of  \cite{Lorentz-1963-Math-Annalen}, we define
$\Delta_n=\Delta_n(t)=\max \{{1}/{n}, \delta_n(t) \}$,
$\delta_n=\delta_n(t)=\sqrt{ {t(1-t)}/{n}},$
and $T_{ns}(x)=\sum_{k=0}^n (k-nx)^s p_{nk}(x)$, $s=0,1,\ldots.$
It is clear that  $\Delta_n(t)=n^{-1}$  if $n\le 4$. For $n>4$,
if $|t-0.5|\le 0.5\sqrt{1-{4}/{n}}$ then $\Delta_n(t)=
    \delta_n(t)$, or $=n^{-1}$.
    If $d=1$ and $f\in C^{(r)}[0,1]$ for $r\ge 2$,  Theorem 1 of \cite{Lorentz-1963-Math-Annalen} gives the  estimate
    $|f(t)-P_n(t)|\le C_r  \Delta_n^{r}(t)\omega_r(\Delta_n(t))$, for $t\in[0,1]$, and $n\ge 1$.
So when $|t-0.5|> 0.5\sqrt{1-{4}/{n}}$  this estimate becomes uniform $|f(t)-P_n(t)|\le C_r  n^{-r}\omega_r(1/n)$.

In order to get a non-uniform estimate,  we need to prove an improved version of Lemma 1 of \cite{Lorentz-1963-Math-Annalen}. It is convenient to denote $\bar T_{ns}(x)= n^{-s}T_{ns}(x)$ and $\bar T_{ns}^*(x) =n^{-s}T_{ns}^*(x):=n^{-s}\sum_{k=0}^n |k-nx|^s p_{nk}(x)$, $s=0,1,\ldots$.
\begin{lemma}\label{lem: estimate of Tns*}
For $s\ge 0$ and some constant $A_s$
\begin{equation}\label{ineq: Tns* improved}
    \bar T_{ns}^*(x)\le A_s  \delta_n^{2\wedge s}(x) \Delta_n^{0\vee(s-2)}(x),
\end{equation}
where $a\vee b=\max(a,b)$, and $a\wedge b=\min(a,b)$. Particularly $A_0=A_1=A_2=1$, $A_3=2$ and $A_4=4$. The equality holds when $s=0,2$.
\end{lemma}
\begin{rem} Lemma 1 of \cite{Lorentz-1963-Math-Annalen} gives
    $\bar T_{ns}^*(x)\le A_s  \Delta_n^s(x)$, $s\ge1$,
which
does not imply zero estimates  at  $x=0,1$.
\end{rem}
\noindent{\bf Proof:} The special results for $s=0,1,2$ are obvious. By the formulas on P. 14 of \cite{Lorentz-1986-book-bernstein-poly} 
we have
$\bar T^*_{n4}(x)=\bar T_{n4}(x)
=n^{-2}\delta_n^2(x)[3n(n-2)\delta_n^2(x)+1]\le 4 \delta_n^2(x)\Delta_n^2(x).$
By the Schwartz inequality
 $\bar T^*_{n3}(x)\le [\bar T^*_{n2}(x) \bar T^*_{n4}(x)]^{1/2}= \delta_n(x)[\bar T^*_{n4}(x)]^{1/2} \le 2 \delta_n^2(x)\Delta_n(x).$
For $s\ge 4$,
Romanovsky \citep[see  Eq.5 of][]{ROMANOVSKY01121923}  has   proved that both $T_{n,2r}(x)$ and $T_{n,2r+1}(x)$
can be expressed as
$nx(1-x)\sum_{l=0}^{r-1}[nx(1-x)]^lQ_{rl}(x),$  where $Q_{rl}(x)$ are   polynomials in $x$  with
coefficients depending on $r$ and $l$  only. Similar to \cite{Lorentz-1963-Math-Annalen}, this implies that
$\bar T_{n,2r}^*(x)=\bar T_{n,2r}(x)\le
A_{2r} \delta_n^2(x) \Delta_n^{2r-2}(x).$
By Schwartz inequality again
$\bar T^*_{n,2r+1}(x) \le [\bar T_{n2}(x)\bar T_{n,4r}(x)]^{1/2}$ $\le
A_{2r+1} \delta_n^2(x)\Delta_n^{2r-1}(x).$
The proof of the Lemma is complete.

Now, Theorem 1 of \cite{Lorentz-1963-Math-Annalen} can be enhanced and generalized as follows.
\begin{theorem}\label{thm: another generalization of thm 1 of Lorentz 1963}
(i) If   $f\in C^{(r)}[0,1]^d$, $r=0,1$,  then
\begin{equation}\label{eq2: approx when r=0,1}
    |f(\bm x)-B^f_{\bm n}(\bm x)| \le 
    (d+1)\omega^{(r)}\left(\max_{1\le j\le d}\delta_{n_j}(x_j)\right) \Big[\sum_{j=1}^d \delta_{n_j}(x_j)\Big]^r,\quad 0\le \bm x\le 1.
\end{equation}
(ii) Let $ r\ge 0$, $m>0$, $M_{\bm i}\ge 0$, be given. Then there exists a constant $C_{r,d}=C_{r,d}(m, \mathscr M_{r})$ such that for each
  function $f(\bm x)\in \Lambda_{r}^{(d)}(m,  \mathscr M_{r})$
one can find a sequence $P_{\bm n}(\bm x)$, $\bm n\ge 1$, of polynomials  with positive coefficients of degree $\bm n$ satisfying
\begin{equation}\label{eq3: approx of poly w pos coeff}
    |f(\bm x)-P_{\bm n}(\bm x)|\le
    C_{r,d} \omega^{(r)}[D_{\bm n}(\bm x)]D_{\bm n}^{r-2}(\bm x)  \Big[ \sum_{j=1}^d   \delta_{n_j}(x_j)\Big]^2,\quad 0\le \bm x\le 1,
\end{equation}
where $D_{\bm n}(\bm x)=\mathop{\max}_{1\le j\le d}\Delta_{n_j}(x_j)$.
\end{theorem}
\begin{rem} Estimates in (\ref{eq2: approx when r=0,1}) are generalizations of (6) and (7) of \cite{Lorentz-1963-Math-Annalen}.
    If $d=1$ and $r\ge 2$, then (\ref{eq3: approx of poly w pos coeff}) is an improved version of Theorem 1 of \cite{Lorentz-1963-Math-Annalen}:
    \begin{equation}\label{eq: improved approx of poly w pos coeff}
    |f(t)-P_n(t)|\le C_r \delta_n^{2}(t)\Delta_n^{r-2}(t)\omega_r(\Delta_n(t)),\quad 0\le t\le 1,\quad n=1,\ldots.
\end{equation}
This indicates that the the approximation  $P_n$ for $f$ performs especially good at the boundaries because the errors are zero at $t=0,1$. However,  results of \cite{Lorentz-1963-Math-Annalen} do not imply this when $r\ge 2$.
\end{rem}
\noindent\textbf{Proof:}
Similar to \cite{Lorentz-1963-Math-Annalen}, we want to prove that, for $r\ge 0$, there exist  polynomials of the form
\begin{equation}\label{eq2: Qnr(f)}
    Q_{\bm n  r}^f(\bm x)=\sum_{\bm k=\bm 0}^{\bm n}\biggr\{f(\sfrac({\bm k},{\bm n}))\!+\!\sum_{i=2}^r \frac{1}{i!}\sum_{\langle\bm i\rangle= i}\!\!{\langle\bm i\rangle\choose \bm i}f^{(\bm i)}(\sfrac({\bm k},{\bm n}))\!\!\prod_{j=1}^d \frac{1}{n_j^{i_j}}\tau_{ri_j}(x_j, n_j) \biggr\}p_{\bm n,\bm k}(\bm x),
\end{equation}
where ${\langle \bm i\rangle \choose \bm i}={\langle \bm i\rangle \choose i_1,\ldots,i_d}$ is the multinomial coefficient, and  $\tau_{ri}(x, n)$'s are polynomials, independent of $f$, in $x$ of degree $i$, in $n$ of degree $\lfloor i/2\rfloor$,
such that for each function $f\in C^{(r)}[0,1]^d$,
\begin{equation}\label{ineq2: for |f-Qnr(f)|}
    |f(\bm x)-Q_{\bm n r}^f(\bm x)|\le C'_{r,d}    \omega^{(r)}[D_{\bm n}(\bm x)] D_{\bm n}^{0\vee(r-2)}(\bm x)  \Big[ \sum_{j=1}^d   \delta_{n_j}(x_j)\Big]^{2\wedge r}
\end{equation}
with $C'_{r,d}$ depending only on $r$ and $d$. 

If $f\in C^{(r)}[0,1]^d$, $r\ge1$,  by the Taylor expansion of $f(\bm k/\bm n)$ at $\bm x$, we have
$$f(\bm x)= f(\sfrac(\bm k,\bm n))-\sum_{i=1}^r \frac{1}{i!}\sum_{\langle\bm i\rangle= i}{\langle\bm i\rangle\choose \bm i}\prod_{j=1}^d(\sfrac({k_j},{n_j})-x_j)^{i_j}
f^{(\bm i)}(\bm x)$$
$$\hspace{7em}+\frac{1}{r!}\biggr\{\sum_{\langle\bm i\rangle= r}{r\choose \bm i}\prod_{j=1}^d(\sfrac({k_j},{n_j})-x_j)^{i_j}\big[f^{(\bm i)}(\bm x) - f^{(\bm i)}(\bm \xi_{\bm k}^{(r)}) \big]\biggr\},$$
where $\bm\xi_{\bm k}^{(r)}$ is on the  line segment connecting $\bm x$ and $\bm k/\bm n$. This equation is also true when $r=0$ by defining $\bm\xi_{\bm k}^{(0)}=\bm k/\bm n$  and the empty sum to be zero. Multiplying both sides by $p_{\bm n,\bm k}(\bm x)$ 
and taking summation over  $\bm 0\le \bm k\le \bm n$, we obtain
\begin{equation}\label{Taylor expansion of f}
f(\bm x)=B^f(\bm x)-\sum_{i=2}^r \frac{1}{i!}\sum_{\langle\bm i\rangle= i}{\langle\bm i\rangle\choose \bm i}\prod_{j=1}^d\bar T_{n_ji_j}(x_j)f^{(\bm i)}(\bm x)+R_{\bm n}^{(r)}(\bm x),
\end{equation}
where 
$r\ge0$, empty sum is zero, and
$$R_{\bm n}^{(r)}=\frac{1}{r!}\biggr\{\sum_{\langle\bm i\rangle= r}{r \choose \bm i}\sum_{\bm k=0}^{\bm n}\prod_{j=1}^d\frac{1}{n_j^{i_j}}\left( k_j-n_jx_j\right)^{i_j}p_{n_jk_j}(x_j)\big[f^{(\bm i)}(\bm x) - f^{(\bm i)}(\bm \xi_{\bm k}^{(r)})\big]\biggr\}.$$
For each $\delta>0$, define
$\lambda=\lambda(\bm x,\bm y;\delta)=\left\lfloor  {|\bm x-\bm y|}/{\delta}\right\rfloor$, where $\lfloor x\rfloor$ is the integer part of $x\ge 0$.
Then $\lambda \delta\le |\bm x-\bm y|<(\lambda+1)\delta$, and
 for   $g\in C[0,1]^d$,
 $|g(\bm x)-g(\bm y)|\le 
(\lambda+1) \omega(g,\delta).$

If  $f\in C^{(r)}[0,1]^d$, $r=0,1$,  then similar to the proofs of Theorems 1.6.1 and 1.6.2 of \cite[][pp. 20-- 21]{Lorentz-1986-book-bernstein-poly} and by (\ref{Taylor expansion of f}) we have
$|f(\bm x)-B^f_{\bm n}(\bm x)|
=|R_n^{(r)}(\bm x)|$.
Because $\lambda(\bm x,  \bm k/\bm n;\delta)
\le \delta^{-1}\sum_{j=1}^d |k_j-n_j x_j|/n_j$, by Lemma \ref{lem: estimate of Tns*} with $s=0,1,2$, we have
\begin{align}\nonumber
  |f(\bm x)-B^f_{\bm n}(\bm x)|
&\le \sum_{\langle\bm i\rangle= r}\omega_{\bm i}(\delta)\biggr[\prod_{j=1}^d \bar T^*_{n_j,i_j}(x_j)
\!\!+\!\!\frac{1}{\delta}  \sum_{l=1}^d \bar T^*_{n_l,i_l+1}(x_l)\!\!
\!\!\mathop{\prod_{1\le j\le d}}_{j\ne l}\!\!\bar T^*_{n_j,i_j}(x_j) \biggr]\\\label{ineq r=0,1}
&\le \sum_{\langle\bm i\rangle= r}\omega_{\bm i}(\delta)\biggr[\prod_{j=1}^d\delta_{n_j}^{i_j}(x_j)
 +\frac{1}{\delta}  \sum_{l=1}^d \delta_{n_l}^{i_l+1}(x_l)\!\!\mathop{\prod_{1\le j\le d}}_{j\ne l}\delta_{n_j}^{i_j}(x_j) \biggr].
\end{align}
The estimates in (\ref{eq2: approx when r=0,1}) follow from (\ref{ineq r=0,1})  with $\delta=\max_{1\le j\le d}\delta_{n_j}(x_j)$. This also proves  (\ref{ineq2: for |f-Qnr(f)|}) with $r=0,1$ and $Q_{\bm n r}^f=B^f$.

If $r\ge 2$, then  we have
\begin{align*}
 |R_{\bm n}^{(r)}|
  &\le  \frac{1}{r!}\biggr\{\sum_{\langle\bm i\rangle= r}{r \choose \bm i}\omega_{\bm i}(\delta)\biggr[\prod_{j=1}^d\bar T^*_{n_ji_j}(x_j)
  +\frac{1}{\delta}\sum_{l=1}^d  \bar T^*_{n_l,i_l+1}(x_l)\!\!\mathop{\prod_{1\le j\le d}}_{j\ne l}  \bar T^*_{n_ji_j}(x_j)\biggr]  \biggr\} 
\\
 \quad\quad &\le  \frac{1}{r!}\biggr\{\sum_{\langle\bm i\rangle= r}{r \choose \bm i}\omega_{\bm i}(\delta)\biggr[\prod_{j=1}^d A_{i_j} \delta_{n_j}^{2\wedge i_j}(x_j) \Delta_{n_j}^{0\vee (i_j-2)}(x_j)   \\
  &\quad +\frac{1}{\delta}\sum_{l=1}^d  A_{i_l+1}\delta_{n_l}^{2\wedge(i_l+1)}(x_l)\Delta_{n_l}^{0\vee(i_l-1)}(x_l)\!\!\!\mathop{\prod_{1\le j\le d}}_{j\ne l} \!\!\! A_{i_j}\delta_{n_j}^{2\wedge i_j}(x_j) \Delta_{n_j}^{0\vee (i_j-2)}(x_j)\biggr]\biggr\}.
\end{align*}
Choosing 
$\delta= D_{\bm n}(\bm x)$,  we have
\begin{align*}
 |R_{\bm n}^{(r)}| &\le  \omega^{(r)}(\delta) \frac{1}{r!}\biggr\{\sum_{\langle\bm i\rangle= r}{r \choose \bm i} \prod_{j=1}^d A_{i_j} \delta_{n_j}^{2\wedge i_j}(x_j) \Delta_{n_j}^{0\vee (i_j-2)}(x_j)   \\
  &\quad + \sum_{l=1}^d\sum_{\langle\bm i\rangle= r}{r \choose \bm i}   A_{i_l+1}\delta_{n_l}^{2\wedge i_l}(x_l)\Delta_{n_l}^{0\vee(i_l-2)}(x_l)\!\!\!\!\mathop{\prod_{1\le j\le d}}_{j\ne l}\!\!\!\! A_{i_j}\delta_{n_j}^{2\wedge i_j}(x_j) \Delta_{n_j}^{0\vee (i_j-2)}(x_j)\biggr\}\\
  &\le   C(r,d) \omega^{(r)}(\delta)\mathop{\max}_{1\le j\le d} \Delta_{n_j}^{r-2}(x_j)  \Big[ \sum_{j=1}^d   \delta_{n_j}(x_j)\Big]^2.
\end{align*}

Similar to \cite{Lorentz-1963-Math-Annalen} we shall prove the existence of $Q_{\bm n r}^f$ by induction in $r$. Assuming that all $Q_{\bm n i}^f$ for $i<r$
are established, we iteratively define
\begin{equation}\label{eq2: iteration for Qnr(f)}
    Q_{\bm n r}^f(\bm x)= B^f(\bm x)-\sum_{i=2}^r \frac{1}{i!}\sum_{\langle\bm i\rangle= i}{\langle\bm i\rangle\choose \bm i}\prod_{j=1}^d\bar T_{n_ji_j}(x_j)  Q_{\bm n, r-i}^{f^{(\bm i)}}(\bm x).
\end{equation}
By Lemma \ref{lem: estimate of Tns*} and the inductive assumption, (\ref{ineq2: for |f-Qnr(f)|}) is  satisfied by (\ref{eq2: iteration for Qnr(f)}) as following.
\begin{align*}
  |f(\bm x)-Q_{\bm n r}^f(\bm x)|
   &\le  \sum_{i=2}^r \frac{1}{i!}\sum_{\langle\bm i\rangle= i}\!\!{\langle\bm i\rangle\choose \bm i}\!\!\prod_{j=1}^d\bar T^*_{n_ji_j}(x_j) |f^{(\bm i)}(\bm x)\!-\!Q_{\bm n,r-i}^{f^{(\bm i)}}(\bm x)|+|R_{\bm n}^{(r)}(\bm x)| \\
   &\le \sum_{i=2}^r \frac{C''_{r,d}}{i!}D_{\bm n}(\bm x)^{i-2}(\bm x)\biggr[\sum_{j=1}^d\delta_{n_j}(x_j)\biggr]^{2}
\\   &\quad \cdot \omega^{(r)}[D_{\bm n}(\bm x)] D_{\bm n}^{0\vee(r-i-2)}(\bm x)\biggr[ \sum_{j=1}^d   \delta_{n_j}(x_j)\biggr]^{2\wedge(r-i)} +|R_{\bm n}^{(r)}(\bm x)|
    \\
   &\le C'''_{r,d}\omega^{(r)}[D_{\bm n}(\bm x)] D_{\bm n}^{0\vee(r-2)}(\bm x) \biggr[ \sum_{j=1}^d   \delta_{n_j}(x_j)\biggr]^{2}.
\end{align*}

Since $f(\bm x)\ge m>0$, by an obvious generalization of remark (a) on p. 241 of \cite{Lorentz-1963-Math-Annalen} with $\bm h=1/\bm n$ we know that $P_{\bm n+r}(\bm x)=Q_{\bm n r}^f(\bm x)$ is a $d$-variate polynomial of degree $\bm n+ r=(n_1+r,\ldots,n_d+r)$ with positive coefficients
for all $\bm n\ge \bm n_{r}(m, \mathscr M_{r})$ so that $$|f(\bm x)-P_{\bm n+r}(\bm x)|\le C_{r,d}\omega^{(r)}[D_{\bm n}(\bm x)] D_{\bm n}^{r-2}(\bm x) \biggr[ \sum_{j=1}^d   \delta_{n_j}(x_j)\biggr]^{2}. $$
Then (\ref{eq3: approx of poly w pos coeff}) follows for all $\bm n$ and a larger $C_{r,d}$ from  $\Delta_{n_j}=\mathcal{O}(\Delta_{n_j+r})$ for all $r\ge 0$.

\section{Conclusion} We have generalized the univariate polynomials with positive coefficients to the multivariate ones and proved an enhanced generalization of Theorem 1 of G. G. Lorentz \cite{Lorentz-1963-Math-Annalen}. The estimation of the degree of approximation of $f\in C^{(r)}[0,1]^d$ by the polynomials with positive coefficients contains a factor
$[ \sum_{j=1}^d   \delta_{n_j}(x_j) ]^{2} = [ \sum_{j=1}^d  \sqrt{x_j(1-x_j)/{n_j}} ]^{2 \wedge r}$ when $r\ge 1$ which is non-uniform even for $\bm x$ close to the vertices of the unit hypercube $[0,1]^d$.
\section*{References}


\end{document}